\documentclass[preprint,10pt]{elsarticle}

\usepackage{amssymb}
\usepackage{amsmath}
\usepackage{algorithm}
\usepackage{algorithmic}
\usepackage{subfigure}
\usepackage{color,xcolor}
\usepackage{graphicx}
\newtheorem{theorem}{Theorem} 
\newtheorem{remark}{Remark} 
\journal{journal math xxx}

\begin{document}

\begin{frontmatter}



\title{A computational inverse random source problem for elastic waves}


\author[a1]{Hao Gu} 
\author[a1]{Tianjiao Wang} 
\author[a1,b]{Xiang Xu } 
\author[c]{Yue Zhao} 
\affiliation[a1]{organization={School of Mathematical Sciences, Zhejiang University},
            city={Hangzhou},
            country={China}}
\affiliation[b]{organization={Center for Interdisciplinary Applied Mathematics, Zhejiang University},
	city={Hangzhou},
	country={China}}
\affiliation[c]{organization={School of Mathematics and Statistics, and Key Lab NAA-MOE, Central China Normal University},
	city={Wuhan},
	country={China}}
\begin{abstract}

This paper investigates the inverse random source problem for elastic waves in three dimensions, 
where the source is assumed to be driven by an additive white noise.
A novel computational method is proposed for reconstructing the variance of the random source from the correlation boundary measurement of the wave field.
Compared with existing multi-frequency iterative approaches, our method is non-iterative and requires data at only a single frequency. As a result, the computational cost is significantly reduced.
Furthermore, rigorous error analysis is conducted for the proposed method, which gives a quantitative error estimate. Numerical examples are presented to demonstrate effectiveness of the proposed method. Moreover, this method can to be directly applied to stochastic Maxwell equations.

\end{abstract}




\end{frontmatter}



\section{Introduction}
\label{sec1}
The inverse source scattering problem is to determine the unknown source from
measurements of the radiating wave fields. 
It has wide applications in many scientific and industrial areas, such as non-invasive medical imaging, antenna synthesis and geophysical exploration \cite{greenleaf2003selected,fokas2004unique}.
The deterministic inverse source scattering problems have been investigated extensively from both theoretical and numerical aspects \cite{devaney2003nonuniqueness,eller2009acoustic,bao2020stability,entekhabi2020increasing}.
In general, it is known that there is no uniqueness for the inverse source problem at a fixed frequency due to the existence of non-radiating sources \cite{albanese2006inverse,hauer2005uniqueness}. Thus, additional information is required for the source in order to obtain a unique solution, such as to seek the minimum energy solution or delta-type sources \cite{bao2002inverse, marengo2002inverse}.  
From the computational point of view, a more challenging issue is the lack of stability. A small variation of the data might lead to a huge error in the reconstruction. 
To regain the uniqueness and obtain increased stability, 
a commonly used method is to employ multi-frequency data \cite{bao2015inverse}.

Stochastic inverse source problems arise in
situations involving randomness in the studied systems, incomplete knowledge of the system, and large scale span \cite{garnier, Devaney, papan1, papan2, baoChow}.
In these situations, the source is modeled by random processes, such as the Gaussian
random field. Therefore, the governing wave equations are stochastic differential equations. 
Compared with the deterministic cases, they are substantially more difficult due to the randomness. In fact, since the driven source
is a random function, it is less meaningful to find a solution to the inverse source problem for a particular realization of the randomness but to
determine the statistical properties of the random source such as the mean and variance.
Similar to deterministic cases, the use of multi-frequency data is effective in developing numerical reconstruction methods for the inverse random source problems. 
For instance, a regularized Kaczmarz method was developed by adopting multi-frequency scattering data
in \cite{bao2016inverse,bao2017inverse} for acoustic and elastic waves. 

In more recent works, \cite{li2024stabilitywhite, wang2025stability2, wang2025stability} proved uniqueness and increasing stability for the inverse random source scattering problems driven by either additive white noise or generalized microlocally isotropic Gaussian field at a single frequency.
Specifically, the variance of the random source can be uniquely determined by the correlation boundary measurement. 
Motivated by the uniqueness result, and the fact that the access to multi-frequency data in practical applications is usually limited,
we aim to develop an effective numerical method to reconstruct random sources driven 
by white noises for wave scattering problems by only using the single-frequency data. Compared with existing multi-frequency iterative reconstruction methods,
see e.g. \cite{bao2016inverse,bao2017inverse}, our proposed method is non-iterative which thus significantly reduces the computational cost.
We consider the more involved three-dimensional elastic waves to demonstrate the effectiveness of the proposed method in handling complex wave systems in higher dimensions, and expect the developed method 
to be applicable to stochastic inverse source problems for acoustic and electromagnetic waves in three dimensions.

The goal is to determine the variance matrix of the vector-valued white noise source by using the displacement of the random wave field measured
on a boundary enclosing the compactly supported source. Motivated by \cite{wang2025stability}, we consider constructing complex exponential solutions to
the elastic wave equation and deduce integral identities. Then by taking covariance of the integral identities and applying Itô's isometry for white noise,
we establish integral equations which connect the diagonal variance matrix and boundary correlation data. By choosing appropriate complex exponential solutions,
we are able to obtain linear combinations of the Fourier coefficients of entries of the variance matrix, which can be explicitly computed by the boundary correlation data.
Moreover, since the variance matrix is diagonal, we shall establish three such integral equations which form a linear system. By solving the linear system
we may obtain the Fourier coefficients for each entry of the variance matrix. 
Computationally, the exponential solutions shall be carefully selected in order to ensure the stability of the reconstruction. To achieve this goal, we 
employ numerical linear algebra techniques to obtain a well-conditioned linear system.
Moreover,  a regularization scheme is employed by incorporating high-frequency truncation for the Fourier transforms of the variances.
The error estimate of the proposed method is derived. Numerical experiments show that the proposed approach is effective to solve the inverse problem,
which also match the derived error estimate.

The rest of the paper is organized as follows. In Section \ref{model}, we introduce the stochastic elastic wave equation. The reconstruction method is proposed in Section \ref{algo}. Section \ref{ee}  is devoted to the error estimate of our proposed method. Numerical examples are presented in Section \ref{numeric} to illustrate the performance of the proposed method and verify the error estimate.

\section{Model problem}\label{model}
Consider the elastic scattering problem of the three-dimensional stochastic Navier equation in a homogeneous and isotropic medium
\begin{equation}\label{ne}
	\Delta^* \boldsymbol{u}+\kappa^2\boldsymbol{u}=\boldsymbol f \quad\mbox{ in } \mathbb{R}^3,
\end{equation}
where $\Delta^* :=\mu\Delta+(\lambda+\mu)\nabla\nabla\cdot$~, $\boldsymbol{u}:\mathbb{R}^3 \rightarrow \mathbb{C}^3$ is the vector-valued displacement of the random wave field, $\kappa>0$ is the angular frequency, $\lambda$ and $\mu$ are the Lam\'{e} constants satisfying $\mu>0$ and $\lambda+2\mu>0$.
The source $\boldsymbol f= (f_1(x), f_2(x), f_3(x))^{\top}$ is assumed to be a random function driven by an additive white noise and takes the form 
$$\boldsymbol f(x) = \boldsymbol \sigma(x)\dot{\boldsymbol W}_x,$$ 
where 
\[
\boldsymbol \sigma(x) = {\rm diag}(\sigma_1(x), \sigma_2(x),\sigma_3(x))
\]
is a deterministic diagonal matrix function characterizing the strength of the random source
with $\sigma_j \geq 0, j=1,2,3$. We assume that $\boldsymbol \sigma$ is supported in a bounded domain $D$.
Here $\boldsymbol W_x = (W_1(x), W_2(x), W_3(x))^{\top}$ is a three-dimensional three-parameter Brownian sheet, where $W_1(x), W_2(x)$ and $W_3(x)$ are three independent one-dimensional three-parameter Brownian sheets.
$\dot{\boldsymbol W}_x$ is a white noise which can be considered as the derivative of the Brownian sheet $\boldsymbol W_x$ in the sense of distribution. 

Denote $B_R:=\{x\in\mathbb{R}^3 : |x|<R\}$ with $R>0$ being a positive constant  such that $D\subset\subset B_R$.
In the exterior domain $\mathbb{R}^3\setminus \overline{D}$ outside the support of the source,
the radiating field $\boldsymbol u$ can be decomposed into the compressional part $\boldsymbol u_{ p}$ and the shear part $\boldsymbol u_{  s}$ as follows:
$$
\boldsymbol u =\boldsymbol u_{p}+\boldsymbol u_{s},\quad
\boldsymbol {u}_{p}=-\frac{1}{\kappa_{p}^2}\nabla\nabla\cdot\boldsymbol u,\quad
\boldsymbol {u}_{s}=\frac{1}{\kappa_{s}^2}\bf{curl}\,\bf{curl}\boldsymbol u,
$$
where $\kappa_p=\frac{\kappa}{\sqrt{\lambda + 2\mu}}$ is the compressional wavenumber, $\kappa_s=\frac{\kappa}{\sqrt{\mu}}$ is the shear wavenumber.

To ensure the well-posedness of the direct problem, the Kupradze-Sommerfeld radiation condition shall be imposed
\begin{align*}
	\lim_{r\to\infty}{r}(\partial_r\boldsymbol u_{p}-\mathrm{i}\kappa_{p}\boldsymbol u_{ p})=0,\quad
	\lim_{r\to\infty}{r}(\partial_r\boldsymbol u_{s}-\mathrm{i}\kappa_{s}\boldsymbol u_{  s})=0,\quad r=|x|,
\end{align*}
uniformly in all directions $\hat{x}=x/|x|$.
It has been shown in \cite{bao2017inverse} that the direct scattering problem admits a unique mild solution:
\begin{equation}
	\boldsymbol u (x)=\int_{D}\mathbf G(x, y) \boldsymbol \sigma(y){\rm d} {\boldsymbol W}_y, \quad a.s. 
	\label{mild_soution}
\end{equation}
Here, The Green tensor $\mathbf G(x, y)$ has the explicit form
\begin{equation*}
	\mathbf G(x, y)=\frac{1}{\mu}g(x, y; \kappa_{s})\mathbf	I +\frac{1}{\kappa^2}\nabla_{x}\nabla^\top_{x}(g(x, y;\kappa_{s})-g(x,y; \kappa_p)),
\end{equation*}
where $\mathbf I$ is the $3\times 3$ identity matrix, and
\begin{equation*}
	g(x, y; \kappa)=-\frac{1}{4\pi}\frac{e^{\mathrm{i}{\kappa}|x-y|}}{|x-y|}
\end{equation*}
is the fundamental solution for the three-dimensional Helmholtz equation.

In this paper, we consider the inverse problem of determine the variance $\boldsymbol \sigma^2$ for the random source from the measured random wave field on $\partial B_R$ at a single frequency.


\section{Reconstruction method}\label{algo}
In this section, we propose a numerical method to reconstruct the variance of the random source.
We begin by constructing a pair of complex exponential solutions $\boldsymbol U_1$ and $\boldsymbol U_2$ to the Navier equation:
\begin{align}\label{Ul}
\boldsymbol  U_1= {\boldsymbol \eta}_1 e^{\mathrm{i}{\boldsymbol \zeta}_1 \cdot x},\quad \boldsymbol  U_2= {\boldsymbol \eta}_2 e^{\mathrm{i}{\boldsymbol \zeta}_2 \cdot x},
\end{align}
where ${\boldsymbol \eta}_l \cdot {\boldsymbol \zeta}_l =0$ and ${\boldsymbol \zeta}_l \cdot {\boldsymbol \zeta}_l =\kappa^2_s$ for $l=1,2$.
It can be verified that the constructed complex exponential solutions satisfy
 \[
 \Delta^*\boldsymbol U_l + \kappa_s^2 \boldsymbol U_l=0.
 \]
Noting that $\boldsymbol f= \boldsymbol \sigma\dot{\boldsymbol W}_x$ and the three components of $\boldsymbol f$ are mutually independent, we obtain the following identity by Ito isometry: 
\begin{align}
	\label{ito}
	\mathbb E\Big[\int_{\mathbb R^3}\boldsymbol f\cdot \boldsymbol U_1{\rm d}x\int_{\mathbb R^3}\boldsymbol f\cdot \boldsymbol U_2{\rm d}x\Big]
	=\int_{\mathbb R^3} \boldsymbol U_1^{\top} {\rm diag}(\sigma_1^2, \sigma_2^2, \sigma_3^2)\boldsymbol U_2{\rm d}x.
\end{align}

In what follows, we show that the left-hand side of  the above identity can be computed from the correlation boundary data. Specifically, by multiplying both sides of the equation
\eqref{ne} by $\boldsymbol U_l$ and using integrating by parts over $B_R$, we obtain
\begin{align*}
	\int_{\mathbb R^3}\boldsymbol f\cdot \boldsymbol U_l{\rm d}x
	&=
	\int_{B_R}\boldsymbol f\cdot \boldsymbol U_l{\rm d} x \notag\\
	&= \int_{B_R} \left(\Delta^* \boldsymbol{u}+\kappa^2\boldsymbol{u}\right) \cdot \boldsymbol U_l{\rm d} x-\int_{B_R} \left(\Delta^* \boldsymbol{U_l}+\kappa^2\boldsymbol{U_l}\right) \cdot \boldsymbol u{\rm d} x \notag\\
	&=\int_{ B_R}\left(\Delta^* \boldsymbol{u}\cdot \boldsymbol U_l
	-\Delta^* \boldsymbol{U_l}\cdot \boldsymbol u
	\right){\rm d}s(x) \notag\\
	&=\int_{ \partial B_R}\left(
	D\boldsymbol u\cdot \boldsymbol U_l-D\boldsymbol U_l \cdot \boldsymbol u
	\right){\rm	d}s(x),\quad l=1,2,
\end{align*}
where the boundary operator $D\boldsymbol{u}=\mu\partial_{\boldsymbol\nu}\boldsymbol{u}
+(\lambda+\mu)(\nabla\cdot\boldsymbol{u})\boldsymbol\nu$. Hence, we have
\begin{align}\label{correlation}
&\mathbb E\Big[\int_{\mathbb R^3}\boldsymbol f\cdot \boldsymbol U_1{\rm d}x\int_{\mathbb R^3}\boldsymbol f\cdot \boldsymbol U_2{\rm d}x\Big] \notag\\
&=\int_{\partial B_R}\int_{\partial B_R}\mathbb E\Big[(D\boldsymbol u(x)\cdot \boldsymbol U_1(x)
- D\boldsymbol U_1(x)\cdot \boldsymbol u(x))\notag\\
&\quad\times(D\boldsymbol u(y)\cdot \boldsymbol U_2(y)
- D\boldsymbol U_2(y)\cdot \boldsymbol u(y))\Big]{\rm d}s(x){\rm d}s(y)\notag\\
&=\int_{\partial B_R}\int_{\partial B_R}\mathbb E\Big[(D\boldsymbol u(x)\cdot \boldsymbol U_1(x))(D\boldsymbol u(y)\cdot \boldsymbol U_2(y))\Big]{\rm d}s(x){\rm d}s(y)\notag\\
&\quad-\int_{\partial B_R}\int_{\partial B_R}\mathbb E\Big[(D\boldsymbol u(x)\cdot \boldsymbol U_1(x))(D\boldsymbol U_2(y)\cdot \boldsymbol u(y))\Big]{\rm d}s(x){\rm d}s(y)\notag\\
&\quad-\int_{\partial B_R}\int_{\partial B_R}\mathbb E\Big[(D\boldsymbol U_1(x)\cdot \boldsymbol u(x))(D\boldsymbol u(y)\cdot \boldsymbol U_2(y))\Big]{\rm d}s(x){\rm d}s(y)\notag\\
&\quad+\int_{\partial B_R}\int_{\partial B_R}\mathbb E\Big[(D\boldsymbol U_1(x)\cdot \boldsymbol u(x))(D\boldsymbol U_2(y)\cdot \boldsymbol u(y))\Big]{\rm d}s(x){\rm d}s(y)\notag\\
&=:\text{correlation data},
\end{align}
where we have denoted the sum of the four integrals by correlation data.

Now we extract the Fourier transforms of the variances $\sigma_j^2$ from the integral $\int_{\mathbb R^3} \boldsymbol U_1^{\top} {\rm diag}(\sigma_1^2, \sigma_2^2, \sigma_3^2)\boldsymbol U_2{\rm d}x$ by choosing appropriate complex exponential solutions. To do this, 
denote $\boldsymbol\eta_1=(\eta_{11}, \eta_{12}, \eta_{13})^\top$ and $\boldsymbol \eta_2=(\eta_{21}, \eta_{22}, \eta_{23})^\top$. Then noting the explicit forms \eqref{Ul} of $\boldsymbol U_l$ we have
\begin{align}\label{IEE}
	&\int_{\mathbb R^3} \boldsymbol U_1^{\top} {\rm diag}(\sigma_1^2, \sigma_2^2, \sigma_3^2)\boldsymbol U_2{\rm d}x \notag\\
	=&\int_{\mathbb{R}^3}\left(\eta_{11} \eta_{21} \sigma_1^2+\eta_{12} \eta_{22} \sigma_2^2+\eta_{13} \eta_{23} \sigma_3^2\right) e^{\mathrm{i}\left(\boldsymbol\zeta_1+\boldsymbol\zeta_2\right)\cdot x} {\rm d} x \notag\\
	=&\eta_{11} \eta_{21} \widehat{\sigma_1^2}(-\boldsymbol\xi)+\eta_{12} \eta_{22} \widehat{\sigma_2^2}(-\boldsymbol\xi)+\left.\eta_{13} \eta_{23} \widehat{\sigma_3^2}(-\boldsymbol\xi)\right|_{\boldsymbol\xi=\boldsymbol\zeta_1+\boldsymbol\zeta_2},
\end{align}
which is indeed a linear combination of Fourier transform of the variance at point $-\boldsymbol\xi=-(\boldsymbol\zeta_1+\boldsymbol\zeta_2)$.
Here $\widehat v$ stands for the Fourier transform of the function $v$ defined by $\widehat v(\boldsymbol\xi)=\int_{\mathbb R^3} v(x) e^{-{\rm i}x \cdot \boldsymbol\xi}\,\mathrm{d}x$. Combing \eqref{ito}--\eqref{IEE}, we establish the following integral  equation
which connects the Fourier transform of entries of the covariance matrix and
the correlation boundary data:
\begin{align}\label{IE}
\eta_{11} \eta_{21} \widehat{\sigma_1^2}(-\boldsymbol\xi)+\eta_{12} \eta_{22} \widehat{\sigma_2^2}(-\boldsymbol\xi)+\left.\eta_{13} \eta_{23} \widehat{\sigma_3^2}(-\boldsymbol\xi)\right|_{\boldsymbol\xi=\boldsymbol\zeta_1+\boldsymbol\zeta_2}=\text{correlation data}.
\end{align}

Notice that in order to further obtain the Fourier transform of each $\widehat{\sigma_j^2}, j=1, 2, 3$, we need to choose three pairs of complex exponential solutions 
$\boldsymbol  U_1^{(k)}= {\boldsymbol \eta}_1^{(k)} e^{\mathrm{i}{\boldsymbol \zeta}_1^{(k)} \cdot x}$, $\boldsymbol  U_2^{(k)}= {\boldsymbol \eta}_2^{(k)} e^{\mathrm{i}{\boldsymbol \zeta}_2^{(k)} \cdot x}$ satisfying ${\boldsymbol \zeta}_1^{(k)}+{\boldsymbol \zeta}_2^{(k)}=\boldsymbol{\xi}, k=1,2,3$,
which yield three equations and form a linear system as follows
\begin{align*}
	\begin{bmatrix}
		\eta^{(1)}_{11} \eta^{(1)}_{21} & \eta_{12}^{(1)} \eta_{22}^{(1)} & \eta_{13}^{(1)} \eta_{23}^{(1)}\\[5pt]
		\eta^{(2)}_{11} \eta^{(2)}_{21} & \eta_{12}^{(2)} \eta_{22}^{(2)} & \eta_{13}^{(2)} \eta_{23}^{(2)}\\[5pt]
		\eta^{(3)}_{11} \eta^{(3)}_{21} & \eta_{12}^{(3)} \eta_{22}^{(3)} & \eta_{13}^{(3)} \eta_{23}^{(3)}\\
	\end{bmatrix}
	\begin{bmatrix}
		\widehat{\sigma_1^2}\\[5pt]
		\widehat{\sigma_2^2}\\[5pt]
		\widehat{\sigma_3^2}
	\end{bmatrix}
	=\begin{bmatrix}
		\text{correlation $\rm data^{(1)}$}\\[5pt]
		\text{correlation $\rm data^{(2)}$}\\[5pt]
		\text{correlation $\rm data^{(3)}$}
	\end{bmatrix}.
\end{align*}
We call the above $3\times 3$ matrix the coefficient matrix and denote it by $A$.
The Fourier transform of $\widehat{\sigma_j^2}, j=1, 2, 3$, can be obtained by solving the linear system.

It should be noticed that even though there are infinitely many theoretical choices for 
$\{\boldsymbol  U_1^{(k)}, \boldsymbol  U_2^{(k)}\}$, in practice, we cannot select them arbitrarily, as the error of the correlation data might be amplified tremendously due to large condition number of the coefficient matrix.
This issue could be mitigated by selecting multiple sets and applying the least squares method, which would however greatly increase computational cost.
Based on the above discussion, we should find three sets which will yield a well-conditioned linear system. Our strategy is to derive a real-valued coefficient matrix
whose diagonal entries are maximized. To achieve this, we consider the case where ${\boldsymbol \zeta}_l^{(k)}$ and $\boldsymbol{\eta}_l^{(k)}, l=1,2$ are real vectors.
Specifically, for any given $\boldsymbol\xi$ satisfying $|\boldsymbol \xi| < 2 \kappa_s$ and $\boldsymbol{\xi} \neq \boldsymbol{0}$, noticing $\boldsymbol\xi=\boldsymbol\zeta_1^{(k)}+\boldsymbol\zeta_2^{(k)}$ and $|{\boldsymbol \zeta}_1^{(k)}| = |{\boldsymbol \zeta}_2^{(k)} |=\kappa_s$, the vectors
$\boldsymbol\zeta_1^{(k)}$ and $\boldsymbol\zeta_2^{(k)}$ should take the following form:
\begin{align} \label{zeta}
\boldsymbol\zeta_1^{(k)}=\frac{\boldsymbol\xi}{2}+\sqrt{(\kappa_{s}^2-\frac{|\boldsymbol\xi|^2}{4})}\boldsymbol\alpha^{(k)},\quad
\boldsymbol	\zeta_2^{(k)}=\frac{\boldsymbol\xi}{2}-\sqrt{(\kappa_{s}^2-\frac{|\boldsymbol\xi|^2}{4})} \boldsymbol\alpha^{(k)}, k=1,2,3,
\end{align}
where $\boldsymbol\alpha^{(k)}$ is a vector satisfying $\boldsymbol\alpha^{(k)} \cdot \boldsymbol\xi=0$, $|\boldsymbol\alpha^{(k)}|=1$. In this case, the selection of each $\{\boldsymbol  U_1^{(k)}, \boldsymbol  U_2^{(k)}\}$ pair is equivalent to choosing a set of parameters $\{ \boldsymbol\alpha^{(k)}, \boldsymbol\eta_1^{(k)}, \boldsymbol\eta_2^{(k)} \}$.

\begin{figure}[htbp]
	\centering
	\subfigure[\label{condrandom} Parameters chosen randomly ]{\includegraphics[width=0.45\textwidth]{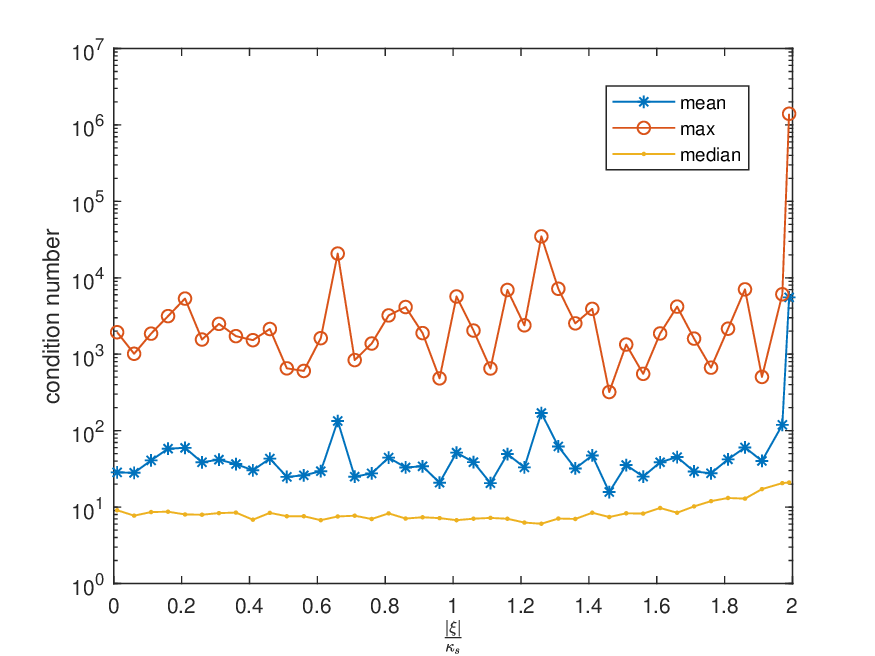}} 
	\subfigure[\label{condp} Parameters chosen by our proposed method ]{\includegraphics[width=0.45\textwidth]{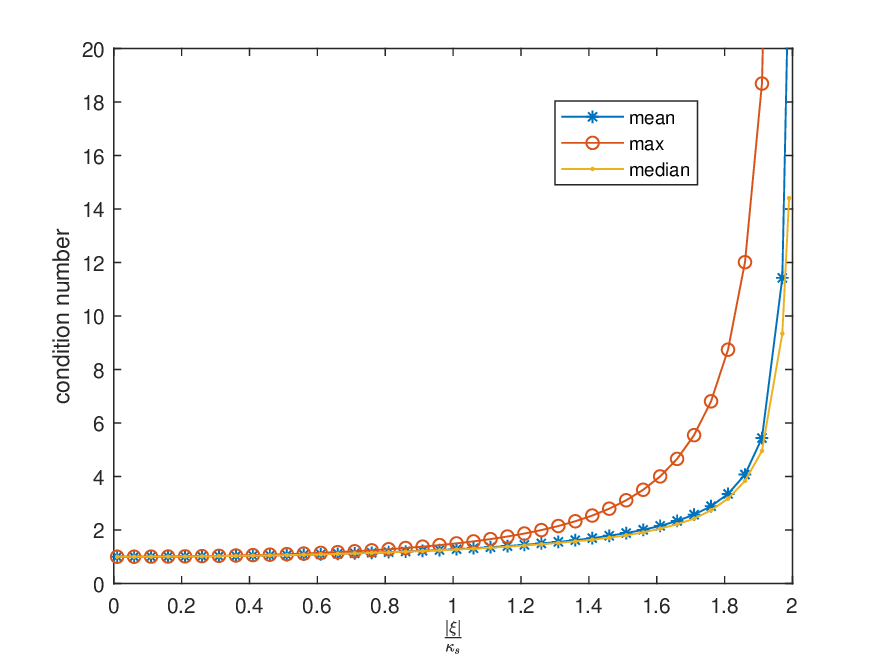}} 
	\caption{Condition number of the coefficient matrices. For each value of $|\boldsymbol\xi|$, we select 256 $\boldsymbol\xi$ uniformly distributed on the sphere with radius $R=|\boldsymbol\xi|$ and determine the corresponding coefficient matrices either randomly or by our proposed approach. For each set of data (each set consists of 256 condition numbers), we compute its mean, maximum, and median.}
	\label{cond}
\end{figure}

Next, we proceed to maximize the diagonal entries of the coefficient matrix. Take the first diagonal element $|\eta^{(1)}_{11}\eta^{(1)}_{21}|$ as an example, since the remainder two diagonal elements could be considered similarly. It is noted that $|\eta^{(1)}_{11}\eta^{(1)}_{21}|$ is determined by the choice of sets $\{ \boldsymbol\alpha^{(1)}, \boldsymbol\eta_1^{(1)}, \boldsymbol\eta_2^{(1)} \}$, or equivalently by the first parameter $\boldsymbol\alpha^{(1)}=(\alpha_1^{(1)},\alpha_2^{(1)},\alpha_3^{(1)})^\top$, because the other parameters $\boldsymbol\zeta_1^{(1)},\boldsymbol\zeta_2^{(1)}$ will also be determined by equation \eqref{zeta}.
Under the conditions $\boldsymbol\eta_l^{(1)}\cdot \boldsymbol\zeta_l^{(1)}=0 ,|\boldsymbol\eta_l^{(1)}|=1, l=1,2$, the maximum values for $|\eta_{11}^{(1)}|$ and $|\eta_{21}^{(1)}|$ are given by
	$$ |\eta_{11}^{(1)}|_{\max}=\frac{\sqrt{\kappa_{s}^2-(\zeta_{11}^{(1)})^2}}{\kappa_{s}}\quad\text{with}
	\quad \zeta_{11}^{(1)}=\frac{\xi_1}{2}+\sqrt{(\kappa_{s}^2-\frac{|\boldsymbol\xi|^2}{4})}\alpha_1^{(1)},
	$$	 
	and
	$$ |\eta_{21}^{(1)}|_{\max}=\frac{\sqrt{\kappa_{s}^2-(\zeta_{21}^{(1)})^2}}{\kappa_{s}}\quad\text{with}
	\quad \zeta_{21}^{(1)}=\frac{\xi_1}{2}-\sqrt{(\kappa_{s}^2-\frac{|\boldsymbol\xi|^2}{4})}\alpha_1^{(1)}.
	$$	
Then we have 
	\begin{align*} &\max |\eta_{11}^{(1)}\eta_{21}^{(1)}| \\ 
	&\quad=\max_{\alpha_1^{(1)}} \frac{\sqrt{\kappa_{s}^2-(\frac{\xi_1}{2}+\sqrt{(\kappa_{s}^2-\frac{|\boldsymbol\xi|^2}{4})}\alpha_1^{(1)})^2}}{\kappa_{s}}\cdot
	\frac{\sqrt{\kappa_{s}^2-(\frac{\xi_1}{2}-\sqrt{(\kappa_{s}^2-\frac{|\boldsymbol\xi|^2}{4})}\alpha_1^{(1)})^2}}{\kappa_{s}}.
	\end{align*}
By direct calculations we have that $|\eta_{11}^{(1)}\eta_{21}^{(1)}|$ reaches a maximum value of $1-\frac{\xi_1^2}{4\kappa_{s}^2}$ when $\alpha^{(1)}_1=0$. 
The other components of $\{ \boldsymbol\alpha^{(1)}, \boldsymbol\eta_1^{(1)}, \boldsymbol\eta_2^{(1)} \}$ can be correspondingly determined.
Specifically, by direct but tedious computations, we obtain
\begin{align*}
	\boldsymbol\alpha^{(1)}=\left(\begin{array}{c}\alpha_1^{(1)}\\ \alpha_2^{(1)}\\\alpha_3 ^{(1)}  \end{array}\right)=\left(\begin{array}{c}0\\ \frac{-\xi_3}{\sqrt{\xi_2^2+\xi_3^2}}\\\frac{\xi_2}{\sqrt{\xi_2^2+\xi_3^2}}  \end{array}\right), \quad \\
	\frac{\boldsymbol\zeta_1^{(1)}}{\kappa_{s}}=\left(\begin{array}{c}\frac{\zeta_{11}^{(1)}}{\kappa_{s}}\\ \frac{\zeta_{12}^{(1)}}{\kappa_{s}}\\\frac{\zeta_{13}^{(1)}}{\kappa_{s}} \end{array}\right)= \left(\begin{array}{c}\frac{\xi_1^{(1)}}{2\kappa_{s}}+\sqrt{1-(\frac{|\boldsymbol\xi|}{2\kappa_{s}})^2}\alpha_1^{(1)} \\ \frac{\xi_2^{(1)}}{2\kappa_{s}}+\sqrt{1-(\frac{|\boldsymbol\xi|}{2\kappa_{s}})^2}\alpha_2^{(1)} \\ \frac{\xi_3^{(1)}}{2\kappa_{s}}+\sqrt{1-(\frac{|\boldsymbol\xi|}{2\kappa_{s}})^2}\alpha_3^{(1)}  \end{array}\right),  \\
	\quad \frac{\boldsymbol\zeta_2^{(1)}}{\kappa_{s}}=\left(\begin{array}{c}\frac{\zeta_{21}^{(1)}}{\kappa_{s}}\\ \frac{\zeta_{22}^{(1)}}{\kappa_{s}}\\\frac{\zeta_{23}^{(1)}}{\kappa_{s}} \end{array}\right)= \left(\begin{array}{c}\frac{\xi_1^{(1)}}{2\kappa_{s}}-\sqrt{1-(\frac{|\boldsymbol\xi|}{2\kappa_{s}})^2}\alpha_1^{(1)} \\ \frac{\xi_2^{(1)}}{2\kappa_{s}}-\sqrt{1-(\frac{|\boldsymbol\xi|}{2\kappa_{s}})^2}\alpha_2^{(1)} \\ \frac{\xi_3^{(1)}}{2\kappa_{s}}-\sqrt{1-(\frac{|\boldsymbol\xi|}{2\kappa_{s}})^2}\alpha_3^{(1)}  \end{array}\right), \\
	\quad \boldsymbol\eta_1^{(1)}=\left(\begin{array}{c}\eta_{11}^{(1)}\\ \eta_{12}^{(1)}\\\eta_{13} ^{(1)}  \end{array}\right)=\frac{1}{{\sqrt{1-(\frac{\zeta_{11}^{(1)}}{\kappa_{s}})^2}}} \left(\begin{array}{c}{1-(\frac{\zeta_{11}^{(1)}}{\kappa_{s}})^2}\\ {\frac{-\zeta_{11}^{(1)}}{\kappa_{s}}\frac{\zeta_{12}^{(1)}}{\kappa_{s}}}\\
		{\frac{-\zeta_{11}^{(1)}}{\kappa_{s}}\frac{\zeta_{13}^{(1)}}{\kappa_{s}}}\end{array}\right), \\
	\quad\boldsymbol\eta_2^{(1)}=\left(\begin{array}{c}\eta_{21}^{(1)}\\ \eta_{22}^{(1)}\\\eta_{23} ^{(1)}   \end{array}\right)=\frac{1}{{\sqrt{1-(\frac{\zeta_{21}^{(1)}}{\kappa_{s}})^2}}} \left(\begin{array}{c}{1-(\frac{\zeta_{21}^{(1)}}{\kappa_{s}})^2}\\ {\frac{-\zeta_{21}^{(1)}}{\kappa_{s}}\frac{\zeta_{22}^{(1)}}{\kappa_{s}}}\\
		{\frac{-\zeta_{21}^{(1)}}{\kappa_{s}}\frac{\zeta_{23}^{(1)}}{\kappa_{s}}}\end{array}\right).
\end{align*} 
For the special case when $\xi_2=\xi_3=0$, $\boldsymbol{\alpha}^{(1)}$ can be taken as $(0, \sin\theta,\cos\theta)^\top$, where $\theta \in [0,2\pi]$ is an arbitrary constant.
Following the same strategy, two other sets of $\{ \boldsymbol\alpha^{(k)}, \boldsymbol\eta_1^{(k)}, \boldsymbol\eta_2^{(k)} \}, k=2, 3,$ can also be chosen to maximize the second and third diagonal elements. 
The curve of the condition number of the coefficient matrix versus $|\boldsymbol\xi|$ is plotted
with a fixed $\kappa_{s}$ in Figure \ref{cond}. 
As shown by the figure, the condition number of the coefficient matrix obtained by our proposed parameter selection method is much smaller than those of other selections.
However, it should be noticed that as $|\boldsymbol{\xi}|$ increases, the condition number also increases accordingly, until it explodes when $|\boldsymbol{\xi}|$ approaches a certain threshold. 
This may be attributed to the intrinsic ill-posedness nature of the inverse problem at a fixed frequency. Therefore, a regularization scheme should be utilized.
To this end, a cutoff frequency is set and the frequencies of the Fourier transforms below this cutoff frequency are used to reconstruct the variance through inverse Fourier transform. This cutoff frequency is set to be proportional to $\kappa_{s}$ which is $\beta\kappa_{s}$ with $\beta\in[0.8,1.25]$ based on our empirical evidence.

As shown in the diagram, Algorithm 1 summarizes the reconstruction procedure. In the following section, an error estimate will be provided for this algorithm.
 
\begin{algorithm}
	\label{algorithm1}
	\caption{Reconstruction of variance at a fixed frequency}
	\begin{algorithmic}
		\REQUIRE
		Boundary random data $\{\boldsymbol u(x_i,\omega_j)\}$, $\{D\boldsymbol u(x_i,\omega_j)\}_{i=1,\cdots, N_{ob},j=1,\cdots, N_s}$.
		$\{x_i\}$,$\{\omega_j\}$ represent observation points on $\partial B_R$ and samples of the random source, respectively. $N_{ob},N_s$ represent the number of observation points and samples, respectively.
	    \STATE Set cutoff frequency $|\boldsymbol\xi|_{\rm max}=\beta \kappa_{s}$ (suggested $\beta\in[0.8,1.25]$), and the set of sample points $\Xi$ in the frequency domain. (initializing)  
		\FOR {$\boldsymbol{\xi}$ in $\Xi$ (loop for Fourier coefficients)}
			\STATE Determine three sets of $\{\boldsymbol \alpha^{(k)},\boldsymbol\eta_1^{(k)},\boldsymbol\eta_2^{(k)} \}$ and form the coefficient matrix $A$ with maximized diagonal elements.
			\FOR {$j=1,\cdots, N_{s}$ (loop for Monte Carlo estimation)}
				\STATE Compute $s_j^{(k)}=  [\Sigma_{i=1}^{N_{ob}} D\boldsymbol u(x_i,\omega_j)\cdot \boldsymbol U_1^{(k)}(x_i)-D\boldsymbol U_1^{(k)}(x_i) \cdot \boldsymbol u(x_i,\omega_j)]\times [\Sigma_{i=1}^{N_{ob}} D\boldsymbol u(x_i,\omega_j)\cdot \boldsymbol U_2^{(k)}(x_i)-D\boldsymbol U_2^{(k)}(x_i) \cdot \boldsymbol u(x_i,\omega_j)], k=1,2,3$.
			\ENDFOR
			\STATE Compute $\text{correlation data}^{(k)} =\frac{16\pi^2 R^4}{N_{ob}^2 N_s}\Sigma_{j=1}^{N_s} s_j^{(k)}, k=1,2,3$.
			\STATE Obtain $\widehat{\sigma_1^2},\widehat{\sigma_2^2},\widehat{\sigma_3^2}$ at $-\boldsymbol{\xi}$ by solving the linear system.
		\ENDFOR
		\STATE Calculate $\sigma_1^2,\sigma_2^2,\sigma_3^2$ using Fourier coefficients at $\Xi$ by inverse Fourier transform.
		\RETURN $\boldsymbol{\sigma}={\rm diag}(\sigma_1,\sigma_2,\sigma_3)$.
	\end{algorithmic}
\end{algorithm}

\section{Error estimate for the algorithm}\label{ee}
In this section, an error estimate is provided for our proposed algorithm. The error associated with the proposed method consists of four parts. 
The first part arises from the high-frequency truncation of the Fourier transforms of the variances, which is utilized as a regularization scheme to ensure the stability of the recovery.
The second part comes from solving the linear system, which results in an interdependence between $\widehat{\sigma_1^2},\widehat{\sigma_2^2},\widehat{\sigma_3^2}$ and thereby induces additional errors. 
The third and fourth part come from the computation of the correlation boundary data, which involve using the Monte Carlo method to calculate expectation and the integration product error of each sample caused by noise in the data, respectively. 

In what follows, we consider the four parts of the error and derive error estimates. We drop the superscripts for different sets of complex exponential solutions
for the sake of convenience. Indeed, the arguments are uniform and can be applied to each case.

\subsection{Truncation of the high-frequency part in the Fourier transform}

In this section, we consider the error resulting from the truncation of the high-frequency part for the Fourier transform of the variances.
Assume that $\|\sigma_j^2\|_{H^1(D)} \leq M$, and denote by $\widetilde{\sigma^2_j}$, the variance after truncating the high-frequency components $|\boldsymbol\xi| > |\boldsymbol\xi|_{\rm max}$, where $j=1,2,3$.
Using Plancherel identity, the error $\boldsymbol{e_1}=(\sigma_1^2-\widetilde{\sigma^2_1},\sigma_2^2-\widetilde{\sigma^2_2},\sigma_3^2-\widetilde{\sigma_3^2})^{\top}$ caused by the high-frequency truncation can be controlled by the regularity of $\sigma_j$:
\begin{align}
	\|\boldsymbol{e_1}\|^2_{L^2(D)}=&\|\boldsymbol{\widehat{e_1}}\|^2_{L^2(\mathbb{R}^3)} \nonumber \\
	=&\sum_{j=1}^{3}\int_{|\xi|>|\xi|_{\max}} |\widehat{\sigma_j^2}|^2{\rm d}\xi \nonumber \\
	=&\sum_{j=1}^{3} \int_{|\xi|>|\xi|_{\max}} \frac{1}{1+|\xi|^2}(1+|\xi|^2) |\widehat{\sigma_i^2}|^2{\rm d}\xi \nonumber\\
	\leq& \frac{3}{1+|\xi|_{\max}^2 }M^2  .\label{err1}
\end{align}

\subsection{Error from solving the linear system}
When solving the linear system for the variances $\widehat{\sigma^2_j}$, the condition number of the coefficient matrix may amplify the error 
in the correlation data by a certain factor.
In this subsection, we analyze this amplification factor, which is indeed the condition number of the coefficient matrix.
From the formulation of $A$, we have that each entry of $A$ depends only on components of the vector $\frac{\boldsymbol{\xi}}{\kappa_s}$.
At each point $\boldsymbol\xi$ we denote the error in the correlation data by $\delta{\boldsymbol{b}}$, and the resulting error in the reconstruction by $\delta{\widehat{\boldsymbol{\sigma}^2}}$.
Then $\delta{\widehat{\boldsymbol{\sigma}^2}}$ admits the following estimate:
\begin{equation}
	\label{err2}
	\|\delta{\widehat{\boldsymbol{\sigma}^2}} \|\leq \frac{{\rm cond}(A)}{\| A\|_2} \|\delta \boldsymbol{b}\|\leq C_1 \|\delta \boldsymbol{b}\|.
\end{equation}
Here, $\|\cdot\|$ denotes the 2-norm of the vector, ${\rm cond}(A)$ is the condition number of matrix $A$ in the 2-norm, $\|A\|_2$ is the spectral norm (induced 2-norm) of matrix $A$ and $C_1$ is a constant depending solely on $\frac{|\boldsymbol\xi|}{\kappa_{s}}$. 
Since we set a cutoff frequency parameter $\beta=\frac{|\boldsymbol\xi|_{\rm max}}{\kappa_{s}}$, $C_1$ can be bounded by a function of $\beta$. 
Despite the absence of a closed-form expression, we can deduce that the function increases monotonically and goes to $+\infty$ as $\beta$ approaches $2$, and can be numerically estimated at any point (similar to Figure \ref{condp}).
For example, if $\beta=1$, then $\rm cond(A)\leq 2$, $\| A\|_2 \geq 0.5$, and thus $C_1$ is bounded by $4$.

\subsection{Error for using the Monte Carlo method to calculate expectation}
The error $\delta{\boldsymbol{b}}$ for obtaining the correlation data consists of two parts.
In this section, we consider the first part $\delta{\boldsymbol{b}_1}$ caused by applying the Monte Carlo method. Typically, the root-mean-square error of the Monte Carlo method decays as $C_2\frac{1}{\sqrt{N_s}}$, where $N_s$ is the number of samples and $C_2<\infty$ is a constant depending on the variance of the estimated quantity. 
In order to ensure the convergence of the algorithm, we need to prove that $\mathrm{Var}\Big[\int_{\mathbb R^3}\boldsymbol f\cdot \boldsymbol U_1{\rm d}x\int_{\mathbb R^3}\boldsymbol f\cdot \boldsymbol U_2{\rm d}x\Big]$ is bounded.
Since $W_1(x), W_2(x)$ and $W_3(x)$ are mutually independent, we have
$$
\begin{aligned}
	&\quad \int_{\mathbb R^3}\boldsymbol f\cdot \boldsymbol U_1{\rm d}x\int_{\mathbb R^3}\boldsymbol f\cdot \boldsymbol U_2{\rm d}x \\
	&=\left[\int_{\mathbb R^3} \left(\sigma_1(x)\dot{W_1}(x)\eta_{11}+\sigma_2(x)\dot{W_2}(x)\eta_{12}+\sigma_3(x)\dot{W_3}(x)\eta_{13}\right)e^{\mathrm{i}{\boldsymbol \zeta}_1 \cdot x}{\rm d}x\right]\cdot \\
	&\quad \left[\int_{\mathbb R^3} \left(\sigma_1(x)\dot{W_1}(x)\eta_{21}+\sigma_2(x)\dot{W_2}(x)\eta_{22}+\sigma_3(x)\dot{W_3}(x)\eta_{23}\right)e^{\mathrm{i}{\boldsymbol \zeta}_2 \cdot x}{\rm d}x\right]\\
	&=\eta_{11}\eta_{21}\int_{\mathbb R^3}\sigma_1(x)\dot{W_1}(x) e^{\mathrm{i}{\boldsymbol \zeta}_1 \cdot x}{\rm d}x\int_{\mathbb R^3}\sigma_1(x)\dot{W_1}(x) e^{\mathrm{i}{\boldsymbol \zeta}_2 \cdot x}{\rm d}x\\
	&\quad+\eta_{12}\eta_{22} \int_{\mathbb R^3}\sigma_2(x)\dot{W_2}(x) e^{\mathrm{i}{\boldsymbol \zeta}_1 \cdot x}{\rm d}x\int_{\mathbb R^3}\sigma_2(x)\dot{W_2}(x) e^{\mathrm{i}{\boldsymbol \zeta}_2 \cdot x}{\rm d}x\\
	&\quad+\eta_{13}\eta_{23} \int_{\mathbb R^3}\sigma_3(x)\dot{W_3}(x) e^{\mathrm{i}{\boldsymbol \zeta}_1 \cdot x}{\rm d}x\int_{\mathbb R^3}\sigma_3(x)\dot{W_3}(x) e^{\mathrm{i}{\boldsymbol \zeta}_2 \cdot x}{\rm d}x.
\end{aligned}
$$ 
Denote
$$
I_j({\boldsymbol \zeta}) = \int_{\mathbb{R}^3} \sigma_j(x) \dot{W_j}(x) e^{\mathrm{i}{\boldsymbol \zeta} \cdot x}{\rm d}x,\quad S=I_j({\boldsymbol \zeta}_1)I_j({\boldsymbol \zeta}_2), \quad j=1, 2, 3.
$$
Noting that $|\eta_{ij}|\leq 1,i=1,2,j=1,2,3$, it suffices to prove $\mathrm{Var}[S]$ is bounded.

It is clear that $I_j({\boldsymbol \zeta}_1)$, $I_j({\boldsymbol \zeta}_2)$ are Gaussian random variables with zero mean.
Using Ito isometry gives 
$$\mathbb{E}[I_j({\boldsymbol \zeta}_1)I_j({\boldsymbol \zeta}_2)]=\int_{\mathbb{R}^3} \sigma_j^2(x) e^{\mathrm{i}({\boldsymbol \zeta}_1+{\boldsymbol \zeta}_2 )\cdot x}{\rm d}x=\widehat{\sigma^2_j}(-{\boldsymbol \zeta}_1-{\boldsymbol \zeta}_2 ),$$
\begin{equation}
	\label{EII}
	\mathbb{E}[I_j({\boldsymbol \zeta}_1)\overline{I_j({\boldsymbol \zeta}_2)}]=\int_{\mathbb{R}^3} \sigma_j^2(x) e^{\mathrm{i}({\boldsymbol \zeta}_1-{\boldsymbol \zeta}_2 )\cdot x}{\rm d}x=\widehat{\sigma^2_j}(-{\boldsymbol \zeta}_1+{\boldsymbol \zeta}_2 ).
\end{equation}


Next, by applying the fourth moment formula for zero-mean jointly Gaussian variables we obtain
$$
\begin{aligned}
\mathrm{Var}(S) &= \mathbb{E}[|S|^2] - |\mathbb{E}[S]|^2 \\
&= \mathbb{E}[|I_j({\boldsymbol \zeta}_1)|^2 |I_j({\boldsymbol \zeta}_2)|^2] - |\mathbb{E}[S]|^2 \\
&=\mathbb{E}[I_j({\boldsymbol \zeta}_1) \overline{I_j({\boldsymbol \zeta}_1)}]\mathbb{E}[I_j({\boldsymbol \zeta}_2) \overline{I_j({\boldsymbol \zeta}_2)}]+
|\mathbb{E}[S]|^2 +|\mathbb{E}[I_j({\boldsymbol \zeta}_1) \overline{I_j({\boldsymbol \zeta}_2)}] |^2- |\mathbb{E}[S]|^2 \\
&=|\widehat{\sigma^2_j}(\boldsymbol{0})|^2+|\widehat{\sigma^2_j}(-{\boldsymbol \zeta}_1+{\boldsymbol \zeta}_2 )|^2\\
&\leq 2|\widehat{\sigma^2_j}(\boldsymbol{0})|^2\\
&< \infty.
\end{aligned}
$$
The last inequality is derived from the fact that $\sigma_j^2$ is a compactly supported function.
Thus, we obtain the following estimate for $ \|\delta{\boldsymbol{b}_1}\|$:
\begin{equation}
	\label{err3}
	 \mathbb{E}[\|\delta{\boldsymbol{b}_1}\|]\leq C_2\frac{1}{\sqrt{N_s}},
\end{equation}
where $C_2=\sqrt{2}\|\boldsymbol{\sigma}\|^2_{L^2(D)}$.

\subsection{Integral product error for each Monte Carlo sample}
For each Monte Carlo sample, the computation of the integral product is imprecise due to the noise in the observation data. In this section, we estimate the integral product error denoted by $\delta{\boldsymbol{b}_2}=(\delta b_{21},\delta b_{22},\delta b_{23})^\top$. Suppose that the observation data $\boldsymbol u^{\delta}(x,\omega_j),D\boldsymbol u^{\delta}(x,\omega_j)$ satisfy   
\[
\|\boldsymbol u^{\delta}(x,\omega_j)-\boldsymbol u(x,\omega_j)\|_{L^1(\partial B_R)} \leq\epsilon,\ \|D\boldsymbol u^{\delta}(x,\omega_j)-D\boldsymbol u(x,\omega_j)\|_{L^1(\partial B_R)} \leq\epsilon,
\]
where $ j=1,2,\cdots,N_s.$
For the sake of convenience, we will omit the sample variable $\omega_j$ in the following discussion since the analysis is the same for all samples.
Denote $$I_l=\int_{ \partial B_R}\left(
D\boldsymbol u\cdot \boldsymbol U_l-D\boldsymbol U_l \cdot \boldsymbol u
\right){\rm	d}s(x),$$
and
$$
I_l^{\delta}=\int_{ \partial B_R}\left(
D\boldsymbol u^{\delta}\cdot \boldsymbol U_l-D\boldsymbol U_l \cdot \boldsymbol u^{\delta}
\right){\rm	d}s(x),$$
where $l=1,2.$
Notice that $\delta b_{21},\delta b_{22},\delta b_{23}$ can be estimated in terms of $|I_1^\delta I_2^\delta-I_1I_2|$.
By direct calculations we have \begin{align}
	\mathbb{E}[|I_1^\delta I_2^\delta-I_1I_2|]&=\mathbb{E}[|(I_1+\epsilon_1)(I_2+\epsilon_2)-I_1I_2|] \notag\\
	&=\mathbb{E}[|\epsilon_1I_2+\epsilon_2I_1+\epsilon_1\epsilon_2|] \notag\\
	&\leq |\epsilon_1| \mathbb{E}[|I_2|]+|\epsilon_2| \mathbb{E}[|I_1|]+|\epsilon_1||\epsilon_2|,\label{4.4.1}
\end{align} where $\epsilon_l=|I_l^\delta-I_l|, \,l=1,2.$
Since $\boldsymbol U_l,l=1,2$ are exponential solutions with a fixed wavenumber $\kappa_{s}$, $|\boldsymbol U_l|, |D\boldsymbol U_l|$ are bounded and thus we have 

\begin{align}
	|\epsilon_l|&=|I_l^\delta-I_l| \nonumber \\
	&=|\int_{ \partial B_R}\left(
	(D\boldsymbol u^{\delta}-D\boldsymbol u)\cdot \boldsymbol U_l-D\boldsymbol U_l \cdot (\boldsymbol u^{\delta}-\boldsymbol{u})
	\right){\rm	d}s(x)|\nonumber\\
	&\leq |\boldsymbol U_l|\|\boldsymbol u^{\delta}-\boldsymbol u\|_{L^1(\partial B_R)}+|D\boldsymbol U_l|\|D\boldsymbol u^{\delta}-D\boldsymbol u\|_{L^1(\partial B_R)} \nonumber\\
	&\leq  D_1\epsilon,\quad l=1,2,\label{D1}
\end{align}
where $D_1$ is a constant depending on $\kappa_{s}$.

Now it remains to prove that $\mathbb{E}[|I_l|]$ is bounded.
We consider an equivalent form of $I_l$:
$$I_l=\int_{\mathbb{R}^3}\boldsymbol f\cdot \boldsymbol U_l{\rm d} x =\sum_{j=1}^3 \eta_{lj} I_j(\boldsymbol\zeta_l), l=1,2.
$$
where
\[
I_j(\boldsymbol\zeta_l) = \int_{\mathbb{R}^3}\sigma_j(x) \dot{W_j}(x) e^{\mathrm{i} \boldsymbol{\zeta}_l \cdot x}{\rm d}x.
\]
We have
$$\mathbb{E}[|I_j(\boldsymbol{\zeta}_l)|^2]=\mathbb{E}[I_j({\boldsymbol \zeta}_l)\overline{I_j({\boldsymbol \zeta}_l)}]=\int_{\mathbb{R}^3} \sigma_j^2(x){\rm d}x.$$
Noting $\sigma_j$ has a compact support, by applying Jensen's inequality to the convex function $f(t) = t^2$ for $t = |I_j(\boldsymbol{\zeta}_l)|$, we obtain
$$
(\mathbb{E} [|I_j(\boldsymbol{\zeta}_l)|])^2 \leq E[|I_j(\boldsymbol{\zeta}_l)|^2] =\int_{\mathbb{R}^3} \sigma_j^2(x){\rm d}x < \infty.
$$
Therefore, $\mathbb{E}[|I_j(\boldsymbol{\zeta}_l)|]$ and thus $\mathbb{E}[|I_l|]$ are bounded, that is
\begin{equation}
	\mathbb{E}[|I_l|] <D_2, \quad l=1,2,
	\label{D2}
\end{equation} 
where $D_2$ is a constant depending on $\|\boldsymbol{\sigma}\|_{L^2(D)}$.
Combining \eqref{4.4.1}--\eqref{D2}, we derive the following error estimate for the integral product: 
$$
\begin{aligned}	
	\mathbb{E}[|I_1^\delta I_2^\delta-I_1I_2|]
	&\leq |\epsilon_1| \mathbb{E}[|I_2|]+|\epsilon_2| \mathbb{E}[|I_1|]+|\epsilon_1||\epsilon_2|\\
	&\leq 2D_1D_2\epsilon+ D_1^2\epsilon^2.
\end{aligned}
$$
By applying the above arguments to each component of $\delta{\boldsymbol{b}_2}$ corresponding to different pairs of exponential solutions, we arrive at the following estimate for $ \|\delta{\boldsymbol{b}_2}\|$:
\begin{equation}
	\label{err4}
	\mathbb{E} [\|\delta{\boldsymbol{b}_2}\|]\leq C_3 \epsilon,
\end{equation}
where $C_3$ is a constant depending on $\kappa_{s}$ and $\|\boldsymbol{\sigma}\|_{L^2(D)}$. 


Combing \eqref{err1}--\eqref{err2}, (\ref{err3}),  and (\ref{err4}), and using the Plancherel theorem
\[
\| \boldsymbol \sigma^2_r-\boldsymbol \sigma^2 \|_{L^2(D)} = \| \widehat{\boldsymbol \sigma^2_r}-\widehat{\boldsymbol \sigma^2} \|_{L^2(\mathbb R^3)},
\]
we derive the following error estimate for our algorithm.
\begin{theorem}[Error estimate for Algorithm 1]
	Denote the reconstructed variance by $\boldsymbol \sigma^2_r$. Assume $\| \sigma_j^2\|_{H^1(D)}\leq M, \|\boldsymbol u^{\delta}-\boldsymbol u\|_{L^1(\partial B_R)} \leq\epsilon,\|D\boldsymbol u^{\delta}-D\boldsymbol u\|_{L^1(\partial B_R)} \leq\epsilon$, and the Monte Carlo sample size is $N_s$.
	Let $|\boldsymbol\xi|_{\rm max}=\beta \kappa_{s}$ with $ \beta\in (0,2)$ be the cutoff frequency.
	 Then the following error estimate holds for Algorithm 1: 
	$$\mathbb{E} \big[\| \boldsymbol \sigma^2_r-\boldsymbol \sigma^2 \|_{L^2(D)}\big] \leq\sqrt{\frac{4\pi}{3}}|\boldsymbol\xi|_{\rm max}^{\frac{3}{2}} C_1(\frac{C_2}{\sqrt{N_s}}+C_3\epsilon)+\frac{\sqrt{3}M}{\sqrt{1+|\boldsymbol\xi|^2_{\rm max}}},
	$$
	where $C_1$ depends solely on $\beta$, and $C_2,C_3$ are constants depending on $M, \kappa_{s}$.
\end{theorem}
\begin{remark}
It can be observed that the error is composed of two parts: the computational error in the low-frequency component and the truncation error in the high-frequency component.
The first part of the error increases with $|\boldsymbol\xi|_{\rm max}$ because the computational domain of low-frequency component is expanding, and the factor $C_1$ is also increasing with $|\boldsymbol\xi|_{\rm max}$. The second part of the error obviously decreases as $|\boldsymbol\xi|_{\rm max}$ grows. 
Therefore, the total error initially increases and then decreases as $|\boldsymbol\xi|_{\rm max}$ increases, which implies
by selecting an appropriate cutoff frequency $|\boldsymbol\xi|_{\rm max}$ the total error can be minimized.
\end{remark}

\section{Numerical examples}\label{numeric}

In this section, we present numerical examples to validate the proposed algorithm and verify the error estimate.

For the discretization of white noise in the random source, we introduce a regular grid of points $\{x_t\}_{t=1}^{N_x}$ with grid step $h$ covering the support of $D$ and approximate the white noise by a piecewise function with Gaussian variables:
$$\dot {W_j}(x)= \sum_{t=1}^{N_x} |K_t|^{-1}\chi_t(x)\int_{K_t}{\rm d} W_j(x)=h^{-\frac{3}{2}}\sum_{t=1}^{N_x}\chi_t(x)Z_{tj},j=1,2,3,$$
where $K_t$ represents a cube centered at $x_t$ with side length $h$, $\chi_t$ denotes the characteristic function of $K_t$, and $Z_{tj}$ are independent Gaussian random variables with mean 0 and variance 1.

The synthetic data for the inverse problem is obtained by taking the convolution of the fundamental solution and the source.
To test the stability of the method, noise with a level of $5\%$ is added to the measurement data
\begin{equation*}
	\boldsymbol u^{\delta}(x_i,\omega_j):=\boldsymbol u(x_i,\omega_j)(1+0.05 \mbox{rand}),\,
	i=1,\cdots,N_{ob}, \, j=1,\cdots,N_{s},
\end{equation*}
\begin{equation*}
	\boldsymbol Du^{\delta}(x_i,\omega_j):=\boldsymbol Du(x_i,\omega_j)(1+0.05 \mbox{rand}),\,
	i=1,\cdots,N_{ob}, \, j=1,\cdots,N_{s},
\end{equation*}
where $\rm rand$ is independently and uniformly distributed random numbers in $[-1,1]$. In the following numerical experiments, we use $N_{ob}=2048$ observation points generated by Fibonacci grid located on $\partial B_R$ with $R=2$.
The number of samples is chosen as $N_s=20000$ such that the error caused by applying the Monte Carlo method is of a lower magnitude than the errors in other steps.

\begin{figure}[htbp]
	\centering
	\includegraphics[width=1.0\textwidth]{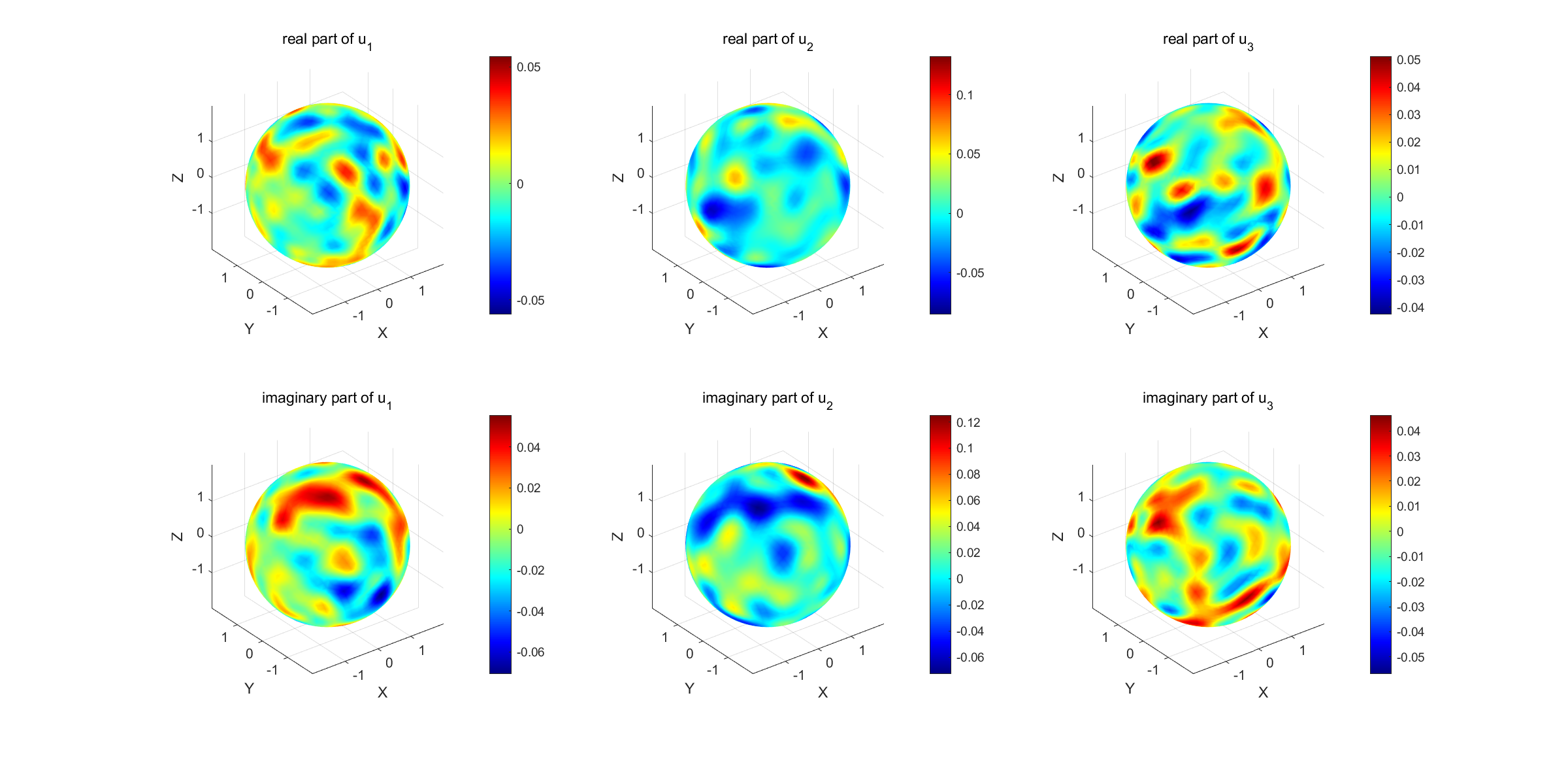}\\
	\caption{Synthetic observation data of a sample.}
	\label{uplot}
\end{figure}

Letting $$\sigma_1=e^{-2(x_1^2+x_2^2+x_3^2)},$$
$$\sigma_2=0.6e^{-8(\sqrt{x_1^2+x_2^2+x_3^2}-0.75)(x_1^2+x_2^2+x_3^2)},$$
$$\sigma_3=0.8e^{-4(x_1^2+(x_2-0.4)^2+(x_3-0.4)^2)}+0.8e^{-4(x_1^2+(x_2+0.4)^2+(x_3+0.4)^2)},$$
the variance $\boldsymbol{\sigma}=\mathrm{diag}(\sigma_1,\sigma_2,\sigma_3)$ is reconstructed inside the domain $D=[-1,1]^3$. The angular frequency $\kappa$ is set as $16$ and the Lamé parameters are taken as $\mu=1,\lambda=2$. Correspondingly, the compressional and shear wavenumbers are $\kappa_p=\kappa/2=8$ and $\kappa_{s}=\kappa=16$, respectively.
We set the discretization grid step of the random source as $h=0.025$ and Figure \ref{uplot} presents the synthetic observation data obtained by solving the forward problem for a sample.

\begin{figure}[htbp]
	\centering
	\subfigure{\includegraphics[width=3.9cm,height=3.2cm]{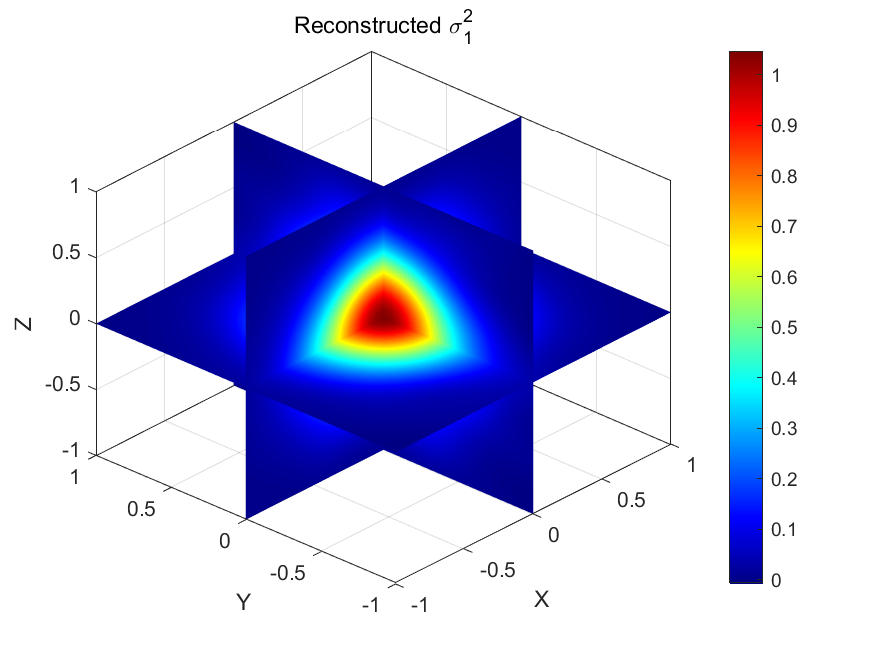}} 
	\subfigure{\includegraphics[width=3.9cm,height=3.2cm]{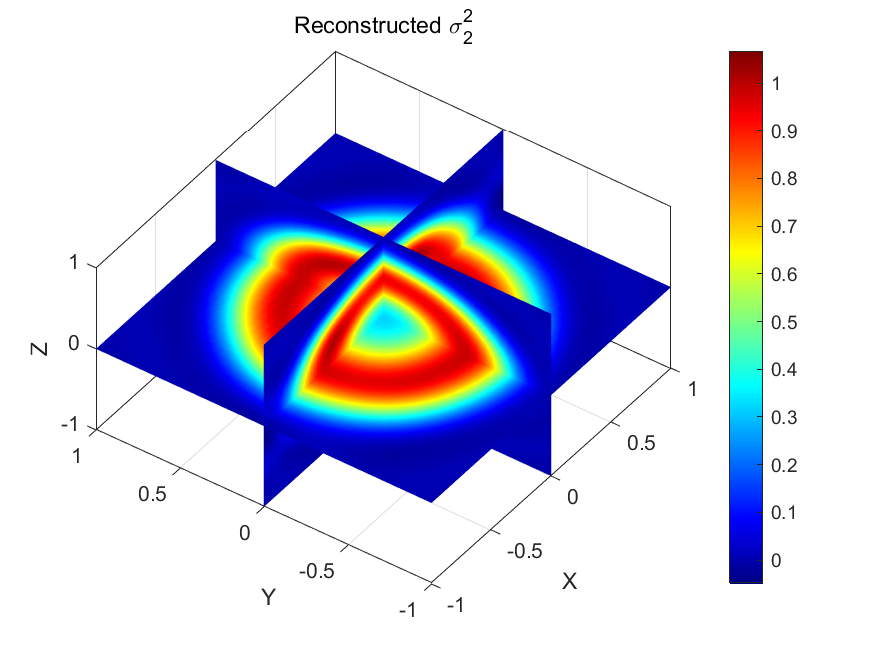}} 
	\subfigure{\includegraphics[width=3.9cm,height=3.2cm]{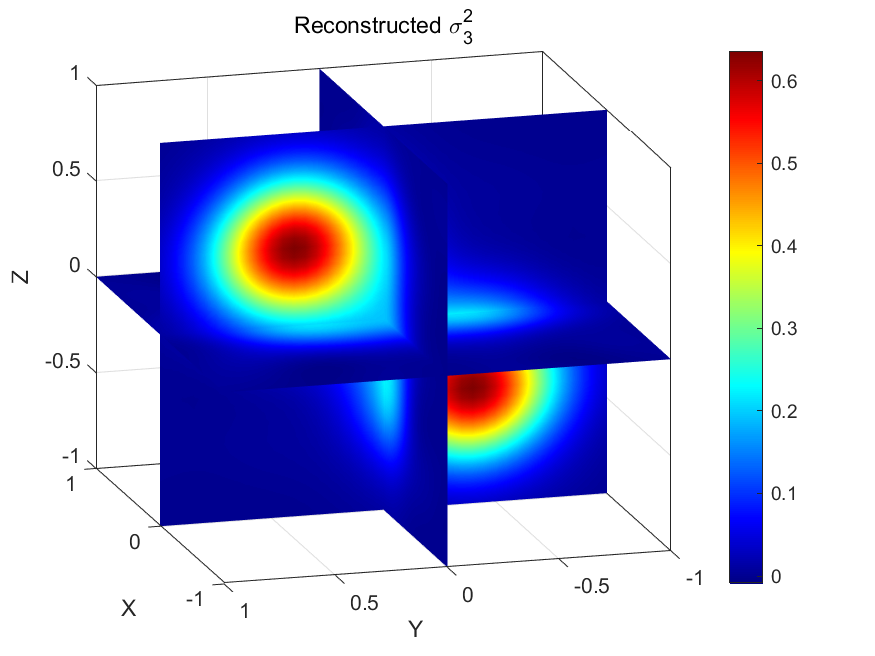}} 
	\caption{Reconstructed variance. From left to right: reconstructed $\sigma_1^2,\sigma_2^2,\sigma_3^2$.}
	\label{reconsturct3d}
\end{figure}

\begin{figure}[htbp]
		\centering
		\includegraphics[width=1.0\textwidth]{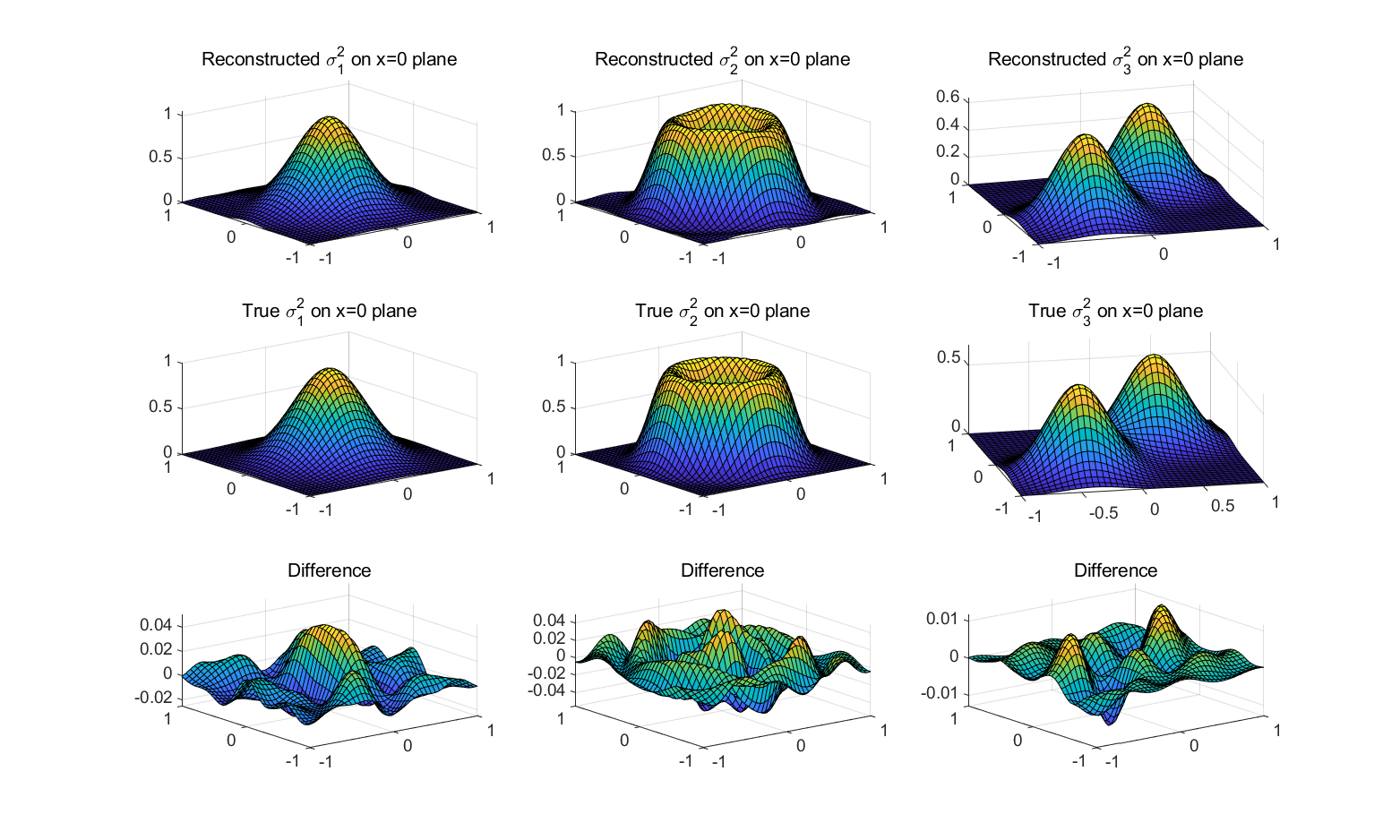}\\
		\caption{Slice of the reconstructed variance on the plane $x=0$. Top row: reconstruction; middle row: ground truth; bottom row: difference.}
		\label{x=0plane}
\end{figure}

\begin{figure}[htbp]
		\centering
		\includegraphics[width=1.0\textwidth]{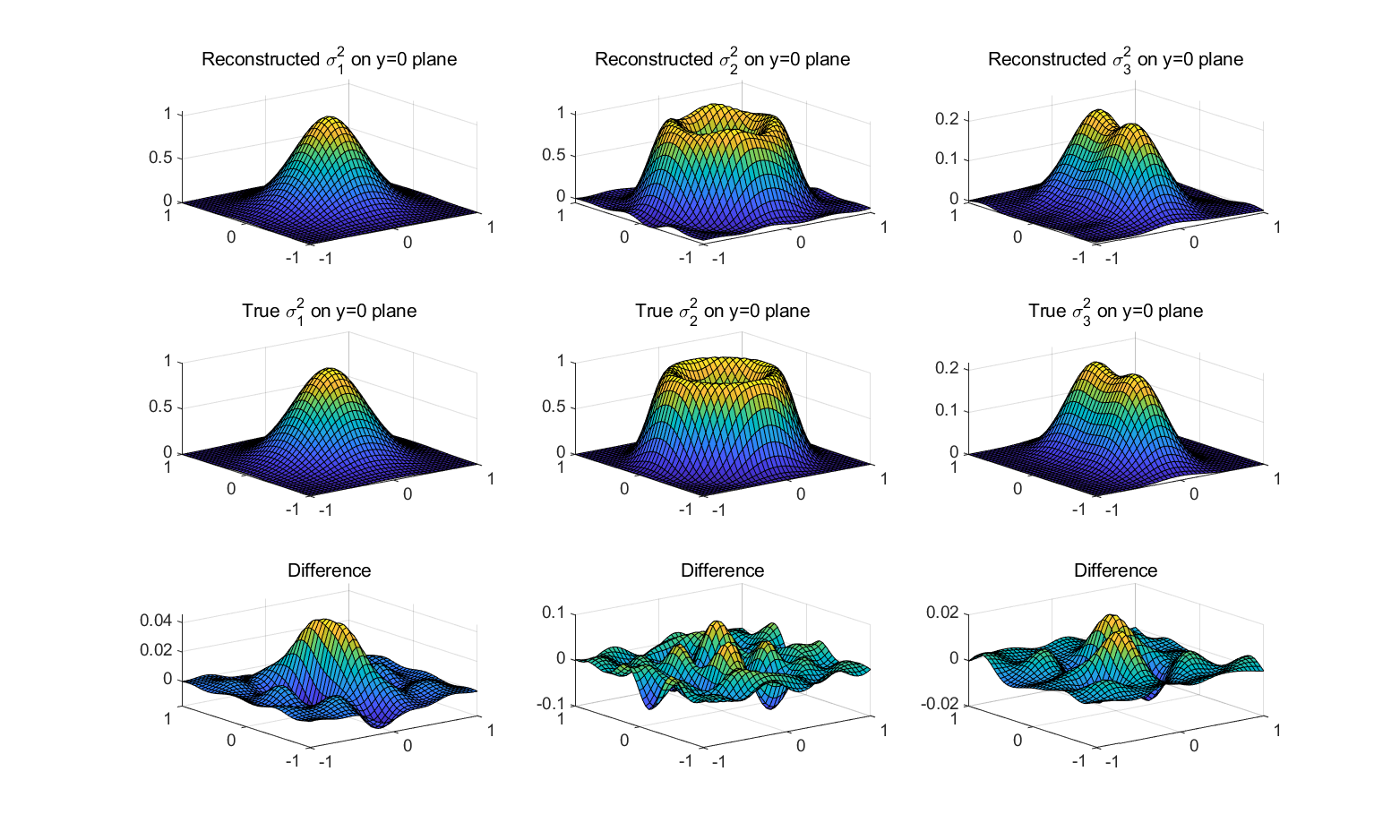}\\
		\caption{Slice of the reconstructed variance on the plane $y=0$. Top row: reconstruction; middle row: ground truth; bottom row: difference.}
		\label{y=0plane}
\end{figure}

\begin{figure}[htbp]
	\centering
	\includegraphics[width=1.0\textwidth]{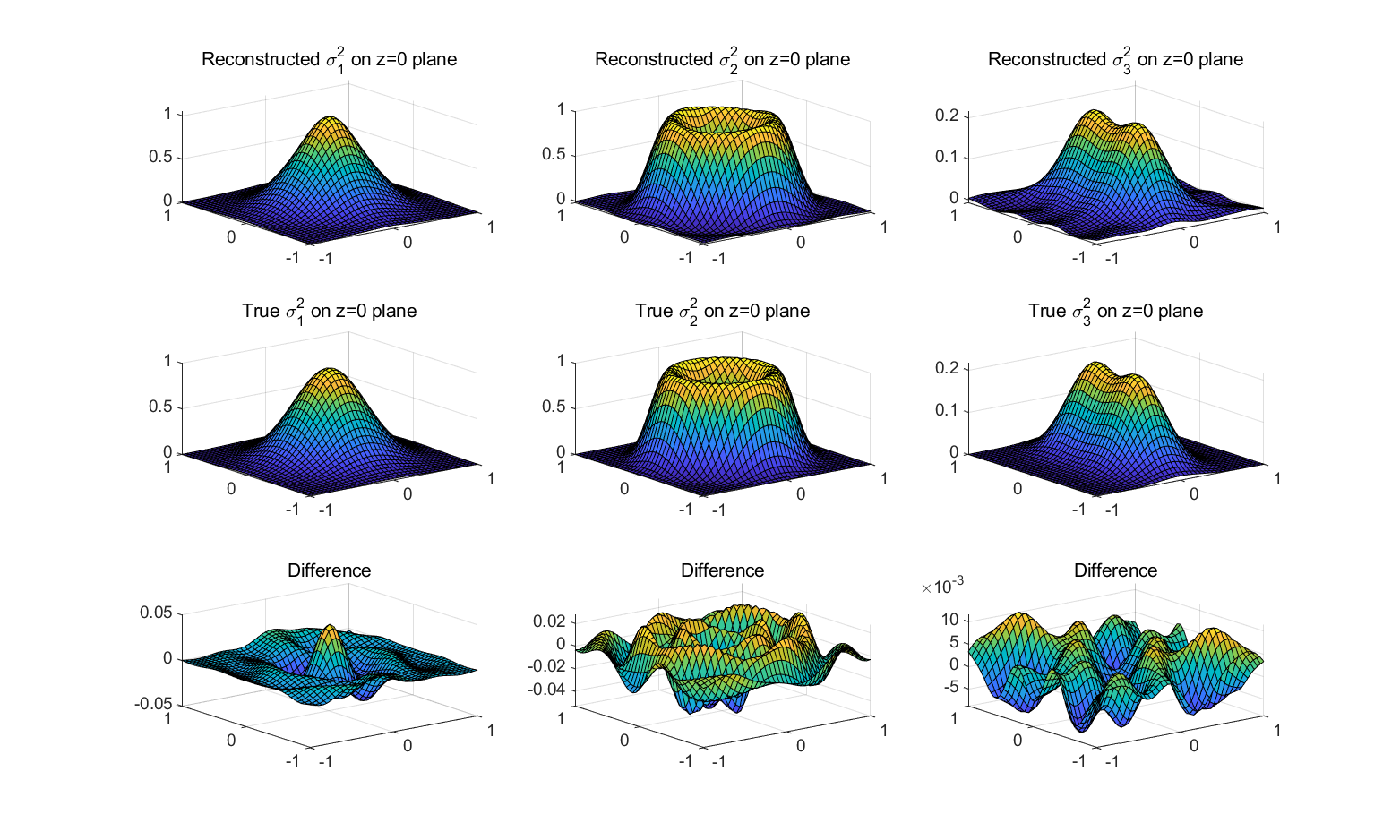}\\
	\caption{Slice of the reconstructed variance on the plane $z=0$. Top row: reconstruction; middle row: ground truth; bottom row: difference.}
	\label{z=0plane}
\end{figure}

We set the truncation frequency threshold $|\boldsymbol\xi|_{\rm max}=14$, i.e., $\beta=0.875,$ and the number of frequency domain samples $|\Xi|\approx 90000$.
Using Algorithm 1, the variance of the random source is reconstructed, and the numerical result is shown in Figure \ref{reconsturct3d}.
We also separately plot the reconstructed values, true values and difference on three slices in Figure \ref{x=0plane}, Figure \ref{y=0plane} and Figure \ref{z=0plane}, repsectively.
Numerical results show that the reconstructed variances not only capture the general outlines, but also match the true values very well in magnitude.
The $L^2(D)$ relative errors of $\sigma_1^2,\sigma_2^2,\sigma_3^2 $ computed by $$err=\frac{\sqrt{\int_{D} (\sigma^2_{r,j}-\sigma_j^2)^2 {\rm d}x }}{\sqrt{\int_{D} (\sigma_j^2)^2 {\rm d}x }},j=1,2,3
$$
are $3.5\%, 4.4\%, 3.9\%$, respectively. The maximum absolute error is less than $0.1$.
	
\begin{figure}[htbp]
	\centering
	\includegraphics[ viewport=85 0 934 285, clip,	width=1.0\textwidth]{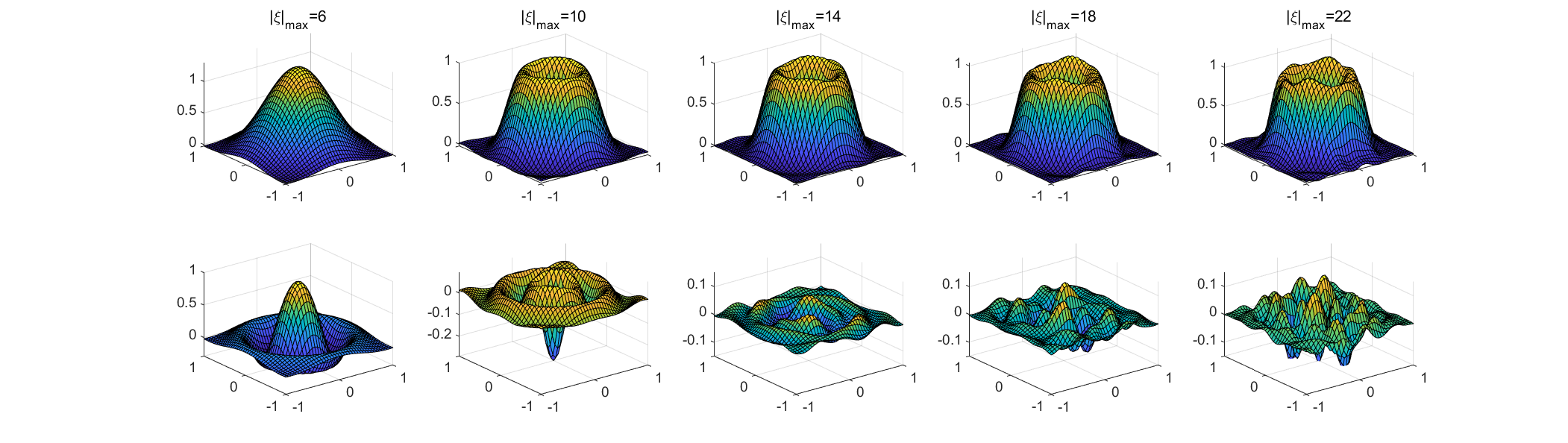}
	\caption{Slice of reconstructed $\sigma_2^2$ on the plane $x=0$ under different $|\boldsymbol{\xi}|_{\rm max}$. Top row: reconstructed results; bottom row: difference.}
	\label{different cutoff frequency parameter }
\end{figure}

The reconstruction method is also tested under different values of $|\boldsymbol\xi|_{\rm max}$. The corresponding errors are shown in Table 1.
\begin{table}[htbp]
	\centering
	\begin{tabular}{|c| c| c| c| c| }
		\hline
		$|\xi|_{\rm max}$ & $\sigma_1^2$ & $\sigma_2^2$ &$\sigma_3^2$ & mean  \\  \hline 
		6 & 13.7\% & 34.9\% & 39.9\% &29.5 \%\\  \hline 
		8 &  2.6\%&  22.4\%& 16.8\%&  13.9\%\\  \hline 
		10 &  1.9\%&  8.4\%& 6.5\%&  5.6\%\\  \hline
		12 &  2.5\%&  7.7\%& 3.4\%&  4.5\% \\  \hline
		\textbf{14} &  \textbf{3.5\%}&  \textbf{4.4\%}& \textbf{3.9\%}&  \textbf{3.9\%} \\  \hline
		16 &  4.8\%&  5.1\%& 5.1\%& 5.0\% \\  \hline
		18 &  6.2\%&  6.1\%& 6.5\%& 6.3\% \\  \hline
		20 &  7.7\%&  7.3\%& 8.0\%&  7.7\% \\  \hline
		22 &  9.5\%&  8.6\%& 10.2\%&  9.4\% \\  \hline
	\end{tabular}
	\caption{ $L^2$ Relative Errors under Different $|\boldsymbol\xi|_{\rm max}$}\label{table1}
\end{table}
It shows that as $|\boldsymbol\xi|_{\rm max}$ increases, the average error initially decreases.
When $|\boldsymbol\xi|_{\rm max}$ reaches $14$, the error decreases to a level close to the minimum.
However, as $|\boldsymbol\xi|_{\rm max}$ continues to increase, the average error begins to rise again.
This matches our theoretical error estimate. In details, as indicated by the error estimate,
the truncation error dominates initially. Then for relatively large $|\boldsymbol\xi|_{\rm max}$,
the error caused by solving the linear system dominates. In our experiment, this part of error attributes to
the error variations of different components of $\boldsymbol{\sigma}$.
Specifically, $\sigma_1$ is the simplest with fewer high-frequency components, while $\sigma_2$ and $\sigma_3$ are more complex, containing more high-frequency components.
Therefore, at low $|\boldsymbol\xi|_{\rm max}$ values, only \(\sigma_1\) can be reconstructed accurately.
As $|\boldsymbol\xi|_{\rm max}$ increases, the condition number of the coefficient matrix $A$ grows, and thus the mutual interference between different components becomes more pronounced.
Therefore, while the reconstruction of $\sigma_2$ and $\sigma_3$ becomes more accurate due to the supplementation of high-frequency components,  the recovery of $\sigma_1$ which lacks high-frequency components is adversely affected. If $|\boldsymbol\xi|_{\rm max}$ is too large, $\sigma_1$, $\sigma_2$ and $\sigma_3$ are all contaminated by amplified noise which leads to increased errors. 
We also plot slices of $\sigma_2^2$ under different values of $|\boldsymbol\xi|_{\rm max}$ in Figure \ref{different cutoff frequency parameter } as a visualization example.
Consistent with previous results, the reconstruction of variances performs best at $|\boldsymbol\xi|_{\rm max}=14$.
When $|\boldsymbol\xi|_{\rm max}=22$, abnormal protrusions can be observed, which are caused by the amplification of noise at high frequencies.
Therefore, we point out that selecting an appropriate cutoff frequency parameter $\beta$ is critical for the reconstruction.

We also test the reconstructed results at different given frequencies. 
For all frequencies, we use the same $\boldsymbol\sigma$, Lamé parameters, Monte Carlo sample size, observation points and noise level as the previous example. 
The same random seeds are used to keep the samples of the random source identical under different frequencies.
The cutoff frequency varies with the given frequency, and we select the optimal cutoff frequency for each one. The corresponding results are presented in Table 2.
It can be observed that the errors decrease as the given frequency increases, which demonstrates the increasing stability. This is consistent with the theoretical results in \cite{wang2025stability}.
\begin{table}[htbp]
	\centering
	\begin{tabular}{|c| c| c| c| c| }
		\hline
		$\kappa$ & $\sigma_1^2$ & $\sigma_2^2$ &$\sigma_3^2$ & mean  \\  \hline 
		2 & 47.5\% & 44.8\% & 68.7\% &53.6 \%\\  \hline 
		4 &  3.9\%&  23.4\%& 18.1\%&  15.1\%\\  \hline 
		6 &  3.4\%&  8.6\%& 7.0\%&  6.3\%\\  \hline
		8 &  4.4\%&  7.9\%& 5.4\%&  5.9\% \\  \hline
		12 &  3.7\%&  4.4\%& 4.7\%& 4.3\% \\  \hline
		16 &  3.5\%& 4.4\%& 3.9\%&  3.9\% \\  \hline

	\end{tabular}
	\caption{ $L^2$ Relative Errors under Different $\kappa$}\label{table2}
\end{table}
\section{Conclusion}
In this paper, we have studied an inverse random source scattering problem for elastic waves, where the source is driven by an additive white noise. 
A stable and efficient numerical method to reconstruct the variance matrix of the random source is proposed, which only employs data at a fixed frequency. 
By constructing pairs of exponential solutions, Fourier coefficients of the variances are related with the correlation data via integral equations.
Computationally, a complete procedure for selecting parameters in the exponential solutions is proposed to achieve stability for the numerical reconstruction.
A quantitative error estimation is derived for the algorithm.
A representative numerical example is reported to demonstrate the reliability and effectiveness of the proposed method, which aligns well with our theoretical error estimate. 

Our algorithm can be directly applicable to Maxwell equations. Consider the time-harmonic Maxwell equations
\begin{align*}
\nabla\times\boldsymbol E(x) - {\rm i}\kappa\boldsymbol H(x) = 0, \quad
\nabla\times\boldsymbol H(x) + {\rm i}\kappa\boldsymbol E(x) = \boldsymbol
f(x),\quad x\in\mathbb R^3,
\end{align*}
where ${\kappa}>0$ is the wavenumber. Eliminating the magnetic field $\boldsymbol
H$ from the Maxwell equations, we obtain the decoupled Maxwell system for the electric
field $\boldsymbol E$:
\begin{equation*}
 \nabla\times(\nabla\times\boldsymbol E)-{\kappa}^2\boldsymbol E={\rm i}{\kappa}\boldsymbol f\quad\text{in}~ \mathbb R^3.
\end{equation*}
Notice that the exponential solution $\boldsymbol  U= {\boldsymbol \eta} e^{\mathrm{i}{\boldsymbol \zeta} \cdot x},$
where ${\boldsymbol \eta} \cdot {\boldsymbol \zeta} =0$ and ${\boldsymbol \zeta}_l \cdot {\boldsymbol \zeta}_l =\kappa^2$, constructed for elastic waves also satisfies the homogeneous
Maxwell equation, that is
\begin{equation*}
 \nabla\times(\nabla\times\boldsymbol U)-{\kappa}^2\boldsymbol U=0\quad\text{in}~ \mathbb R^3.
\end{equation*}
Then repeating the arguments for elastic waves the reconstruction method and the error estimate can be extended to Maxwell equations
in a straightforward way. The extension to the scalar-valued Helmholtz equation also follows the same arguments which is easier.

A possible continuation of this work is to investigate the inverse random source problems in inhomogeneous background medium. In such scenarios, 
explicit exponential solutions are no longer available and new method shall be developed.
We hope to be able to report the progress on these problems in the future.

\section*{Acknowledgement}
The work is partly supported by National Key Research and Development Program of China (No. 2024YFA1012300), National Natural Science Foundation of China (No. 12525112, No. 12371421), the Fundamental Research Funds for the Central Universities, and the Open Research Project of Innovation Center of Yangtze River Delta, Zhejiang University. This work does not have any conflicts of interest.
\bibliographystyle{siam}
\bibliography{refs_abbre}
\end{document}